\documentclass{amsproc}
\usepackage{geometry}                % See geometry.pdf to learn the layout options. There are lots.
\usepackage{graphicx}
\usepackage{amssymb}
\usepackage[final]{pdfpages}
\title[Laplacians on Julia Sets III]{Laplacians on Julia Sets III: Cubic Julia Sets and Formal Matings}

\author[C. Spicer]{Calum Spicer$^1$}
\address{Math Department, Johns Hopkins University}
\email{calum@jhu.edu}
\thanks{$^1$Research supported by the National Science Foundation through the Research Experience for Undergraduates Program at Cornell}

\author[R. S. Strichartz]{Robert S. Strichartz$^2$}
\address{Math Department, Mallot Hall, Cornell University, Ithaca NY 14853}
\email{str@math.cornell.edu}
\thanks{$^2$Research supported in part by the National Science Foundation, grant DMS-0652440}

\author[E. Totari]{Emad Totari$^3$}
\address{Math Department, U.C. Riverside}
\email{etota001@ucr.edu}
\thanks{$^3$Research supported by the National Science Foundation through the Research Experience for Undergraduates Program at Cornell}

\begin{document}
\maketitle
\begin{center}
Dedicated to the memory of Benoit Mandelbrot
\end{center}
\begin{abstract}
We continue the study of constructing invariant Laplacians on Julia sets, and studying properties of their spectra.  In this paper we focus on two types of examples: 1) Julia sets of cubic polynomials $z^3 + c$ with a single critical point; 2) formal matings of quadratic Julia sets.  The general scheme introduced in earlier papers in this series involves realizing the Julia set as a circle with identifications, and attempting to obtain the Laplacian as a renormalized limit of graph Laplacians on graphs derived form the circle with identifications model.  In the case of cubic Julia sets the details follows the pattern established for quadratic Julia sets, but for matings the details are quite challenging , and we have only been completely successful for one example.  Once we have constructed the Laplacian, we are able to use numerical methods to approximate the eigenvalues and eigenfunctions.  One striking observation from the data is that for the cubic Julia sets the multiplicities of all eigenspaces (except for the trivial eigenspace of constants) are even numbers.  Nothing like this is valid for the quadratic julia sets studied earlier.  We are able to explain this, based on the fact that three is an odd number, and more precisely because the dihedral-3 symmetry group has only two distinct one-dimensional irreducible representations.
\end{abstract}

\pagebreak

\section{Introduction}

This is the third of a sequence of papers (following [FS] and [ADS]) exploring the definitions and properties of Laplacians on Julia sets.  One theme is the extension of ideas pioneered by Jun Kigami in the context of post-critically finite (PCF) self-similar fractals ([Ki], [S]).  The Julia sets do not have the strictly self-similar structure, but they do share the finite ramifications properties of PCF fractals.  Kigami's approach is to approximate the fractal by a sequence of graphs, take a sequence of graph energies on these graphs and then an appropriately renormalized limit to obtain an energy on the fractal.  Combining the energy and a measure on the fractal via a weak formulation yields a Laplacian, which can then be described as a suitably renormalized limit of graph Laplacians, echoing the basic idea of calculus that derivatives are limits of difference quotients.  The second theme is that Julia sets may be realized as circles with identifications.  Indeed, the external ray parametrization gives a map from the circle onto the Julia set, and the dynamics on the Julia set is conjugate to a simple map on the circle.  Thus we can hope to take simple structures on the circle and transport them to the Julia set.  For example, the standard Lebesgue measure on the circle may be pushed forward to the Julia set, and this easily provides the measure we will use.  The situation with the energy is more challenging.

We may informally refer to this approach as the {\it method of Peano curves}.  Indeed, the classical Peano curves also give a realization of domains in the plane as circles with identifications.  Of course, no one in their right mind would suggest that the usual planar Laplacian be constructed by pushing forward structures on the circle to planar regions via a Peano curve.  Nevertheless, our approach might be successful in the context of other fractals.

In the first paper [FS] we studied a simple family of quadratic Julia sets, including the familiar Basilica and Douady Rabbit (the Basilica had been treated earlier in [RT]).  In this family there are only a countable number of identifications on the circle, and they arise in a coherent iterative scheme that gives rise to a family of graphs whose vertices are identified points and whose edges are arcs of the circle.   We used a rather ad hoc method to renormalize the graph energies in order to obtain an energy in the limit that is invariant under the doubling map on the circle (conjugate to the dynamics on the Julia set).  In [ADS] the ad hoc procedure was replaced by a more systematic approach and this approach was used to study a few more complicated examples.

We briefly outline the method.  We start with a finite set of identified points on the circle, $X^{(1)}$.  To define an energy on the resulting graph we need to assign conductances (reciprocals of resistances) to the edges.  The simplest choice would be to take the resistance to be equal to the length of the edge, but this is not satisfactory.  Instead, we divide the edges up into different types, and multiply the reciprocal of the length by a factor $b_n$ that depends on the type to obtain the desired conductance.  We then use the dynamics (in the case of quadratic Julia sets this is just the doubling map on the circle) to build a sequence of graphs, $X^{(m)}$ and energies $E_b^{(m)}$ that depend on the choice of weight factors $\{b_n\}$.  The key condition we want is that each level energy restricts to a multiple of the previous level energy, so $E_b^{(m)}(\tilde{u}) = rE_b^{(m-1)}(u)$, when $\tilde{u}$ is the harmonics extension of $u$ (the one that minimizes the $E_b^{(m)}$ energy).  This leads to a nonlinear eigenvalue problem for the $\{b_n\}$.  The key technical problem is to solve this eigenvalue problem.  Note that this is exactly the problem encountered in the case of PCF fractals (see [P] for recent work on that case).

In this paper we study a family of Julia sets for higher order polynomials, but with the condition that the polynomial has a single critical points, with canonical form $z^p+c$.  We are able to generalize the construction in [FS] to this case, and in particular we study in detail some examples when $p = 3$.  Although the construction itself follows the same outline, the structure of the spectrum of the Laplacian is quite different.  One striking feature is that all eigenspaces (beyond the trivial space of constant functions) have even multiplicity.  We will give an explanation for this phenomenon based on the fact that 3 is an odd number, and more precisely the fact that the dihedral-$p$ group for $p$ odd has only two distinct irreducible representation of dimension one (when $p$ is even there are four such representations).

The other example we study is the formal mating of two quadratic Julia sets, specially a dendrite studied in [ADS] and the anti-Rabbit.  We are able to construct the energy explicitly, and we examine in detail the structure on the spectrum of the Laplacian, which is more closely related to the spectrum for the Rabbit than for the Dendrite.  We also briefly discuss two other simple mating for which we were unable to carry out the method outlined above.   This strongly suggests that new ideas will be required to understand these example and other matings.

In section 2 we describe the construction of invariant energies and Laplacians for families of Julia sets for $p-$th order polynomials with a single critical point.  In section 3 we discuss matings of Julia sets.  In section 4 we present numerical data for eigenvalues and eigenfunctions of the Laplacian for these examples.  More data may be found on the website [ http://www.math.cornell.edu/$\sim$etota001/ ].  In section 5 we discuss the structure of the various spectra, giving rigorous explanations of observations that were first made from the experimental results.

\pagebreak

\section{A Simple Family of Julia Sets}

Fix an integer $p\geq 2$, and consider Julia sets of polynomials $z^p+c$.  These polynomials have a single critical points, so their Julia sets are closely related to Julia sets of quadratic polynomials (the case $p = 2$).  Typically the Julia sets may be realized as circles with identifications via exterior ray parameterizations.  Rather than use the parameter $c$ we use a parameter $\theta$ in the circle $\mathbb{R}/\mathbb{Z}$.  The correspondence $c \mapsto \theta$ is not one-to-one, but different choices of $c$ corresponding to the same $\theta$ yield topologically equivalent Julia sets.

Assume for simplicity that $\theta$ is rational.  In fact, we will mainly be dealing with the choice of $\theta = \frac{1}{p(p^k-1)}$ for $k \geq 2$.  Subdivide the circle into the $p$ arcs $A_n = (\theta + \frac{n}{p}, \theta + \frac{n+1}{p}]$.  For each $t$ in the circle define its kneading sequence $(n_0, n_1, ...)$ by the conditions $(p^jt \mod 1) \in A_{n_j}$, $j = 0, 1, ...$, and identify points $t, t'$ if and only if they have the same kneading sequence.  Only rational points get identified, so there are only a countable number of identifications, and this gives a description of the Julia set as a limit of the circle with a finite number of identifications.

For the special case $\theta = \frac{1}{p(p^k-1)}$ we observe that the $k$ points $\frac{p^n}{p^k-1}$, $n = 0, 1, ..., k-1$ are all identified (kneading sequence $(0, 0, 0, ...)$, and these are the only identifications among points of the form $\frac{j}{p^k-1}$.  Among points of the form $\frac{j}{p(p^k-1)}$ we will identify $\{\frac{p^n + \ell (p^k-1)}{p(p^k-1)}\}_{n = 1, ..., k}$ for each $\ell = 0, ..., p-1$ because they have the kneading sequence $(\ell, 0, 0, ...)$.  This produces the first approximating graph $X^{(1)}$ with $p$ vertices.  The edges are the pairs of adjacent vertices on the circle.  Note that $p(p-1)$ of the edges are loops, connecting identified points.  Although the loops do not contribute to the energy, it is convenient to keep them in order to describe the inductive construction of the graph.  Note that the lengths of the edges are $\frac{p^n(p-1)}{p(p^k-1)}$ for $n = 0, ..., k-1$.  We call them edges of type $n$.  For $n \geq 1$ these are loops, while for $n = 0$ these are edges joining $\frac{p^k + \ell (p^k-1)}{p(p^k-1)}$ and $\frac{p + (\ell +1)(p^k-1)}{p(p^k-1)}$.

Inductively, we define the graph $X^{(m)}$ by retaining all the identifications in $X^{(m')}$ for $m' = 1, 2, ..., m-1$ and identifying $\{\frac{p^n + \ell (p^k-1)}{p^m(p^k-1)}\}_{n=1, ..., k}$ for $0 \leq \ell \leq p^m-1$.  These points have kneading sequence $(n_0, n_1, ..., n_{m-1}, 0, 0, ...)$ where $\ell = \sum_{j=0}^{m-1}n_jp^j$.  Moreover, we can give a subdivision rule that describes the edges in $X^{(m)}$ in terms of the edges in $X^{(m-1)}$.   The edges of type $n$ in $X^{(m-1)}$ have length $\frac{p^n(p-1)}{p^{m-1}(p^k-1)}$.  When $n < k-1$ they remain edges in $X^{(m)}$ of type $n+1$.  When $n = k-1$ the edges subdivides into $p-1$ loops of type $n$ for each $n \geq 1$ and $p$ edges of type 0, where the sequence of types is $0, 1, ..., k-1, 0, 1, ..., k-1, ..., 0, 1, k-1, 0$.  Note that the lengths add up correctly:
$$
(p-1)\Big(\frac{p(p-1)}{p^m(p^k-1)}+ \frac{p^2(p-1)}{p^m(p^k-1)}+...+\frac{p^{k-1}(p-1)}{p^m(p^k-1)}\Big)+p\frac{(p-1)}{p^m(p^k-1)} = \frac{p^{k-1}(p-1)}{p^{m-1}(p^k-1)}
$$

Figures 2.1 illustrates the identifications for the cubic Basilica ($p=3, k=2$) and cubic Rabbit ($p=3, k=3$).  Figures 2.2 shows the corresponding Julia sets.

\begin{figure}[h]
\numberwithin{figure}{section}
\caption{Cubic Basilica and Cubic Rabbit Identifications (level 2)}
\centering
\includegraphics[scale = 0.33]{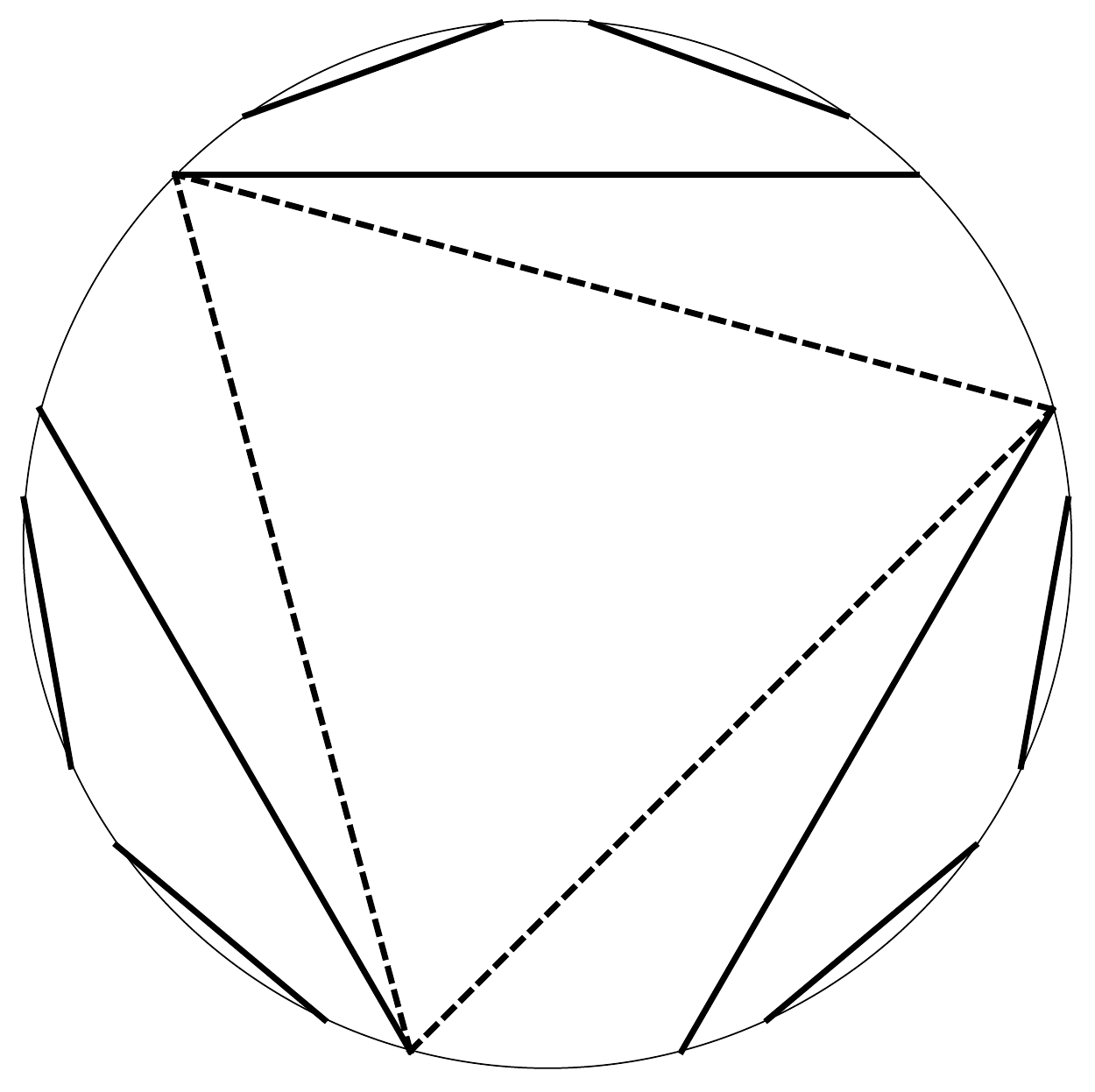} \includegraphics[scale =0.5]{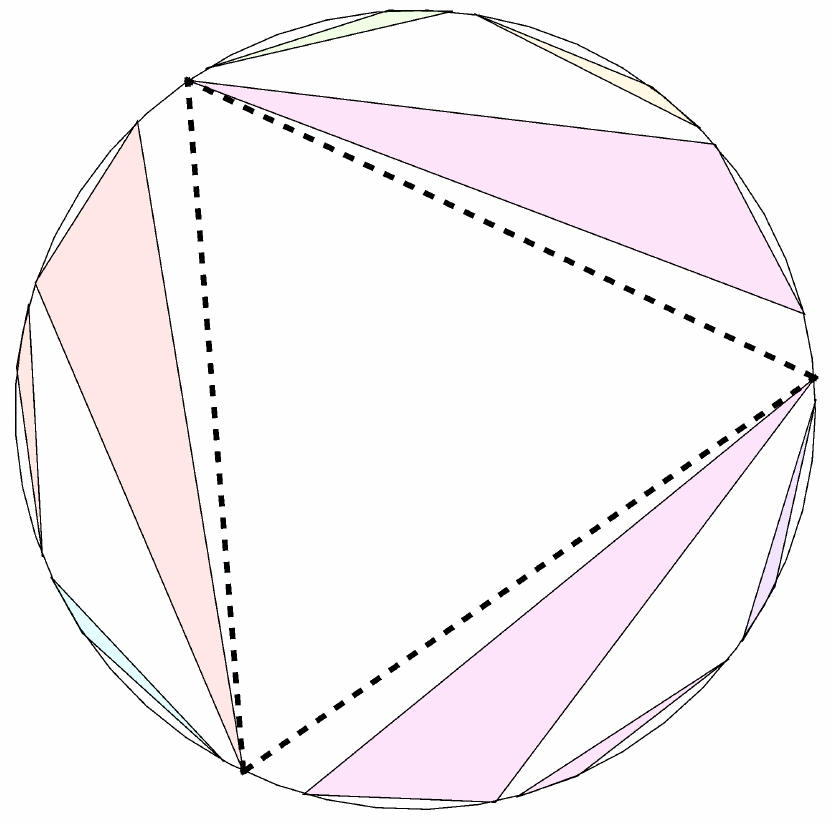}
\end{figure}

\begin{figure}[h]
\numberwithin{figure}{section}
\caption{Cubic Basilica and Cubic Rabbit}
\centering
\includegraphics[scale = 0.5]{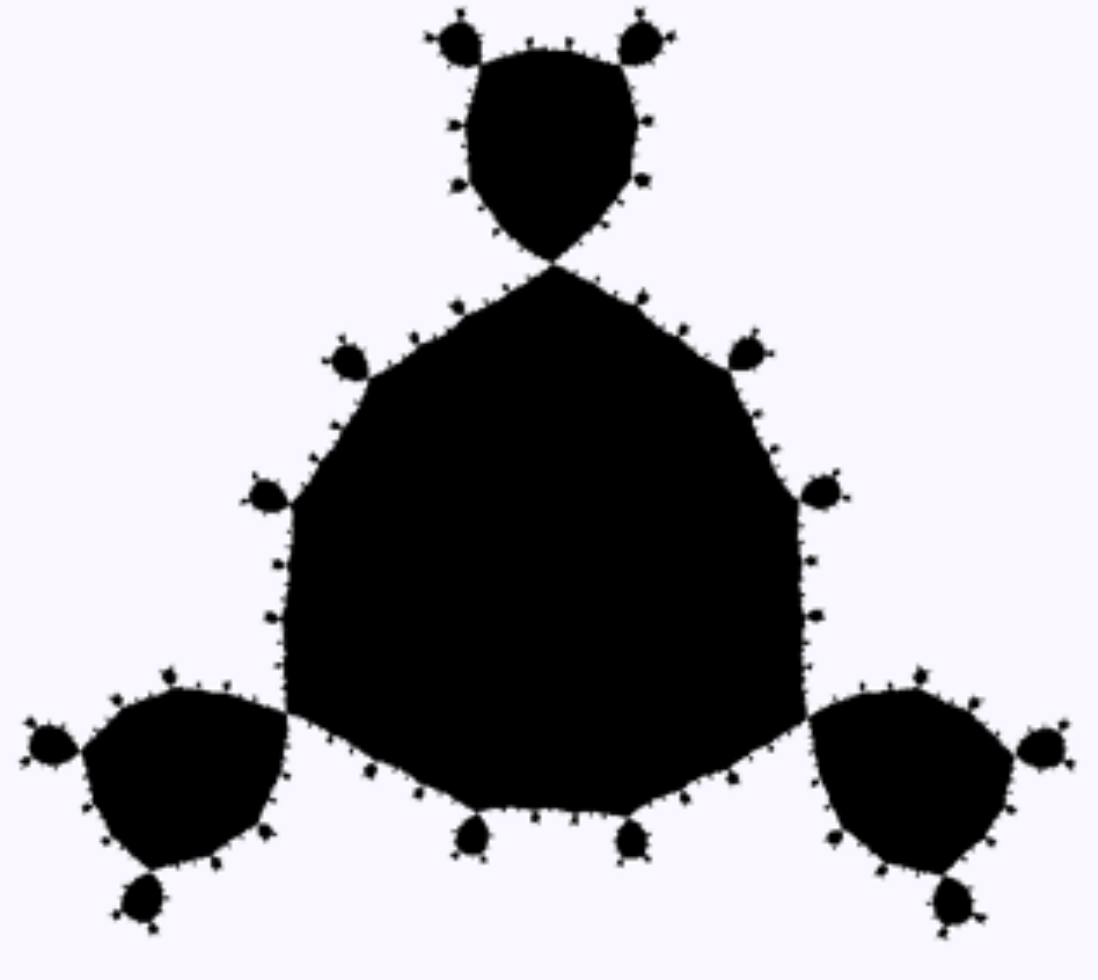} \includegraphics[scale = 0.5]{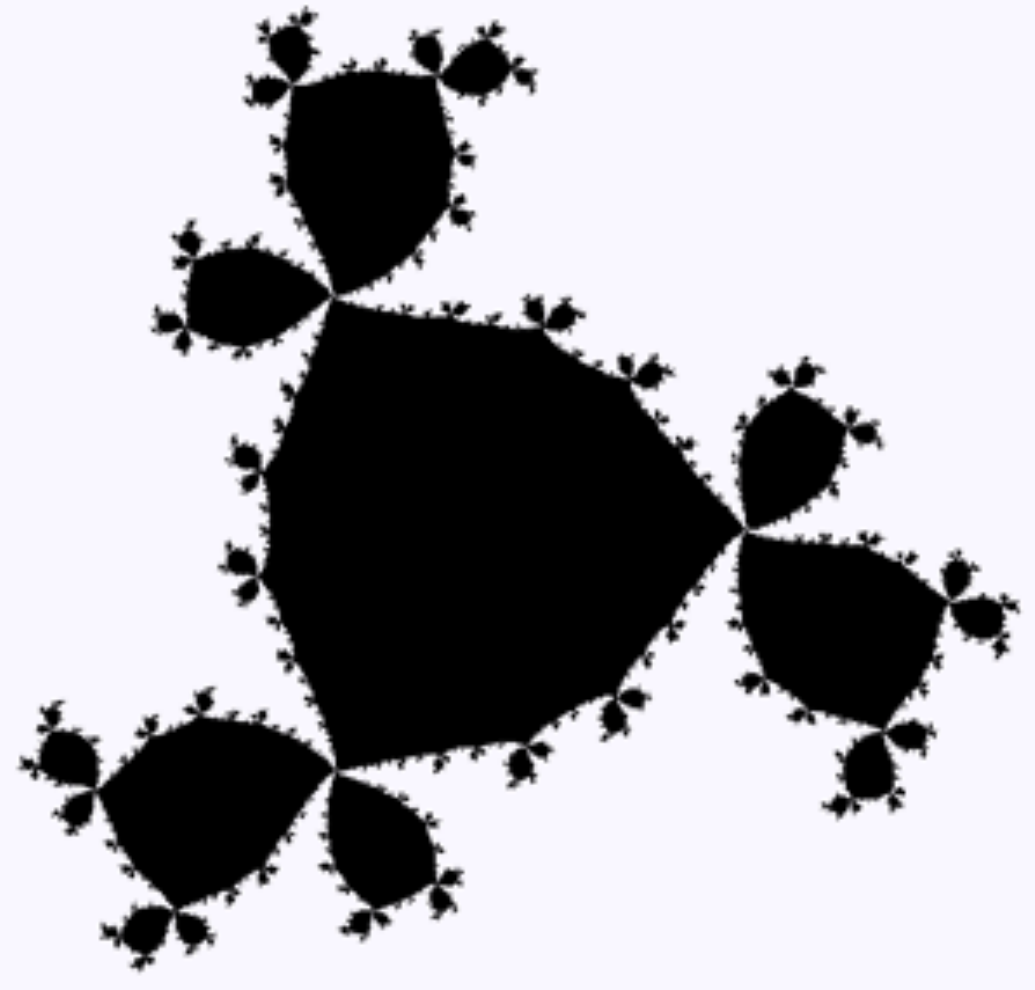}
\end{figure}

We may now describe the energies on the graphs $X^{(m)}$ that we will use to approximate the energy on the Julia set.  Let $u$ be a function on $X^{(m)}$, or equivalently a function on the circle that respects identifications, i.e, $u(x) = u(y)$ if $x$ and $y$ are identified.  Then, let 
\begin{equation*}
E_n^{(m)}(u) = \sum \frac{ |u(x) - u(y)|^2}{ |x - y|}
\tag{2.1}\end{equation*}
summed over all type $n$ intervals $[x, y]$, and
\begin{equation*}
E_b^{(m)}(u) = \sum_{n=0}^{k-1}b_nE_n^{(m)}(u)
\tag{2.2}\end{equation*}
for positive weights $\{b_n\}$.  Note that loops do not contribute to these energies, but we leave them in the sum (2.1) for simplicity of notation.  Also $|x-y| = \frac{p^n(p-1)}{p^m(p^k-1)}$.  If $u$ is defined on $X^{(m-1)}$, we consider all possible extension to $X^{(m)}$ and let $\tilde{u}$ be the one that minimizes $E^{(m)}_b$.  Call this the {\it harmonic extension}.  It might seem that the harmonic extension depends on the choice of weights, but in fact it does not.  In fact is is rather easy to describe the harmonic extension and to relate $E^{(m)}_b(\tilde(u))$ to $E_b^{(m-1)}(u)$.  All the new vertices in $X^{(m)}$ lie in the type $k-1$ edges of $X^{(m-1)}$, and each such edge contains $p-1$ sets of $k$ identified points.  We have to chose the new values of $\tilde{u}$ on these sets of identified points to interpolate linearly between the values of $u$ on the endpoints of the original edge in order to minimize the contribution to $E^{(m)}_0$ from this set of edges.  In other words, if the original edge in $X^{(m-1)}$ was $[x, y]$ with difference $u(y) - u(x)$, we end up with differences $\frac{1}{p}(u(y) - u(x))$ on each of the $p$ edges of type 0.  (The other edges are loops.)  So the edge that contributed $$\frac{|u(y)-u(x)|^2}{(\frac{p^{k-m}(p-1)}{p^k-1})}$$ to $E^{(m-1)}_{k-1}$ now contributes $$p\frac{|\frac{u(y)-u(x)}{p}|^2}{(\frac{p^{-m}(p-1)}{p^k-1})}$$ to $E^{(m)}_0$, a magnification factor of $p^{k-1}$.  There are no other type 0 edges in $X^{(m)}$ and all the type $n$ edges for $n \geq 1$ are just the type $n-1$ edges in $X^{(m-1)}$.  Thus, 
\begin{equation*}
E_{n'}^{(m)}(\tilde{u}) = \sum_{n=0}^{k-1} M_{n'n}E_n^{(m-1)}(u)
\tag{2.3}\end{equation*}
for the matrix
\begin{equation*}
M = \begin{pmatrix}
0 & \cdots & 0 & p^{k-1} \\ 
1 & \ddots & \vdots & 0 \\
\vdots &\ddots& \ddots & \vdots \\
 0 & \cdots & 1 & 0
 \end{pmatrix}
 \tag{2.4}\end{equation*}
By choosing $\{b_n\}$ to be a left eigenvector for $M$, 
\begin{equation*}
rb_n = \sum_{n'=0}^{k-1}b_{n'}M_{n'n}
\tag{2.5}\end{equation*}
we obtain, 
\begin{equation*}
E^{(m)}_b(\tilde{u}) = rE^{(m-1)}_b(u)
\tag{2.6}\end{equation*}
Since $M^k = p^{k-1}I$, it is easy to see that $r = p^{\frac{k-1}{k}}$ and $b = (1, r^{\frac{1}{k}},r^{\frac{2}{k}}, ..., r^{\frac{k-1}{k}})$ give the unique positive solution to (2.5).  So, we make this choice and renormalize by 
\begin{equation*}
\mathcal{E}^{(m)}(u) = r^{-m}E_b^{(m)}(u)
\tag{2.7}\end{equation*}
so that 
\begin{equation*}
\mathcal{E}^{(m)}(\tilde{u}) = \mathcal{E}^{(m-1)}(u).
\tag{2.8}\end{equation*}
Then, 
\begin{equation*}
\mathcal{E}(u) = \lim_{m \rightarrow \infty} \mathcal{E}^{(m)}(u)
\tag{2.9}\end{equation*}
is well-defined on $[0, \infty]$ with $\mathcal{E} = 0$ if and only if $u$ is constant.  We say $u \in dom(\mathcal{E})$ if $\mathcal{E}(u) < \infty$.

Next we consider the behavior of the energy under the action of the map $u \rightarrow u \circ  P$ where $P(x) = px \mod 1$.  Note that $P$ is conjugate to the polynomial $z^p+c$ on the actual Julia set.  Suppose $E$ is an edge of type $n$ in $X^{(m-1)}$.  Then $\frac{1}{p}E + \frac{\ell}{p}$ for $\ell = 0, 1, ..., p-1$ are edges in $X^{(m)}$ also of type $n$.  Since the length of the edges are reduced by the factor $p$, it follows that the contribution from the $p$ edges in $X^{(m)}$ to $E_n^{(m)}(u \circ P)$ is exactly $p^2$ times the contribution from $E$ to $E_n^{(m-1)}(u)$.  Thus, 
\begin{equation*}
E_n^{(m)}(u \circ P) = p^2 E_n^{(m-1)}(u)
\tag{2.10}\end{equation*}
Since (2.10) is independent of $n$ we obtain
\begin{equation*}
E_b^{(m)}(u \circ P) = p^2 E_b^{(m-1)}(u)
\tag{2.11}\end{equation*}
hence,
\begin{equation*}
\mathcal{E}(u \circ P) = \frac{p^2}{r}\mathcal{E}(u) = p^\frac{k+1}{k}\mathcal{E}(u)
\tag{2.12}\end{equation*}

To define a Laplacian we need a measure on the Julia set, which we take to be the Lebesgue measure, $\mu$, on the circle, since the identical points form a set of measure zero.  This is called the {\it equilibrium measure} because
\begin{equation*}
\int f \circ P d\mu = \int f d\mu
\tag{2.13}\end{equation*}

We define $\Delta u$ via the identity
\begin{equation*}
\mathcal{E}(u, v) = - \int (\Delta u)v d\mu
\tag{2.14}\end{equation*}
 for all $v \in dom(\mathcal{E})$, where the bilinear form $\mathcal{E}(u, v)$ is obtained from the quadratic form $\mathcal{E}(u) = \mathcal{E}(u, u)$ and the polarization identity.  It follows from (2.12) and (2.13) that 
  \begin{equation*}
\Delta (u \circ P) = p^{\frac{k+1}{k}}(\Delta u) \circ P
\tag{2.15}\end{equation*}
and this is the invariance condition of $\Delta$ under the action of $P$.
 
We can also represent the Laplacian as a limit of renormalized graph Laplacians $\Delta_m$ on the graphs $X^{(m)}$.  Note that there are exactly $p^m$ distinct vertices in $X^{(m)}$ and they are equally distributed, so the right side of (2.14) is approximated by $-p^{-m}\sum_{x \in X^{(m)}} (\Delta u(x)) v(x)$.  For a fixed $x \in X^{(m)}$, we choose $v = \delta_x$ (in other words $v(x') = \delta_{xx'}$) and define
\begin{equation*}
\Delta_m u(x) = -p^m \mathcal{E}^{(m)}(u, \delta_x)
\tag{2.16}\end{equation*}

We can write this more explicitly by noting that each $x$ has $2k$ edges (some may be counted twice), with two of each $n$ type.  Let $e(x)$ denote this set of edges, and for each edge let $y$ denote the other endpoint (note we may have $y = x$ if the edge is a loop), and let $n$ denote the type.  Then
 \begin{equation*}
\Delta_m u(x) = \sum_{e \in e(x)} c_n (u(y) - u(x))
\tag{2.17}\end{equation*}
with the conductances $c_n$ given by 
 \begin{equation*}
c_n = \frac{p^mr^{-m}b_n}{(\frac{p^n(p-1)}{p^n(p^k-1)})} = p^{m-n+\frac{m+n}{k}}(\frac{p^k-1}{p-1})
\tag{2.18}\end{equation*}

\pagebreak

\section{Formal Matings of Julia Sets}

Suppose $\mathcal{J}'$ and $\mathcal{J}''$ are Julia sets for the polynomials $P(z) = z^p+ c'$ and $P'(z) = z^p+c''$ with the same degree $p$.  If we realize $\mathcal{J}'$ and $\mathcal{J}''$ as circles with identifications corresponding to the parameters $\theta'$ and $\theta''$, then we can realize the formal mating $\mathcal{J}$ of $\mathcal{J}'$ and $\mathcal{J}''$ as a circle with identifications $t_1 \equiv t_2$ if and only if $t_1$ and $t_2$ are identified in $\mathcal{J}'$ or $1-t_1$ and $1-t_2$ are identified in $\mathcal{J}''$.  The main problem we consider is how to construct an invariant (under $t \rightarrow pt \mod 1$) energy on $\mathcal{J}$ if we know how to construct invariant energies on $\mathcal{J}'$ and $\mathcal{J}''$.  More specifically, suppose we have increasing sequences $X'^{(m)}$ and $X''^{(m)}$ of identified points in $\mathcal{J}'$ and $\mathcal{J}''$.  Then we form $X^{(m)} = X'^{(m)} \cup (1- X''^{(m)})$ and make the appropriate identifications.  We will always assume that $X'^{(m)} = p X'^{(m-1)} \mod 1$ and $X''^{(m)} = p X''^{(m-1)} \mod 1$, so that $X^{(m)} = p X^{(m-1)} \mod 1$.  We make $X^{(m)}$ into a graph by declaring edges between all neighboring points on the circle.

Our strategy would be to define different types of edges in these graphs, then form the type $n$ energy via (2.1) and the full energy via (2.2) for a suitable choice of weights.  We would then want to define subdivision rules for how an edges of type $n$ in $X^{(m-1)}$ breaks up into edges of different types in $X^{(m)}$.  We could then find the energy minimizing extension $\tilde{u}$ and compute a matrix $M$ for which (2.3) is valid.  Note that we can expect that the matrix $M$ depends on the weights $\{b_n\}$, so that (2.5) becomes a nonlinear eigenvalue problem.  Nevertheless, if we are able to find a positive solution, then we can define a renormalization (2.7) of the energy (2.2) such that (2.8) holds, and then define an energy on $\mathcal{J}$ via (2.9).

It appears that this strategy is successful in only a few special cases.  We will discuss one example where it works, and then very briefly two examples where we were unable to make it work.  In all these examples $p=2$, so we are considering formal matings of quadratic Julia sets.  From the point of view of dynamical systems, the most interesting formal matings are the ones which yield Julia sets of rational mappings.  These are called simply {\it matings}.  Our main example is not such a mating.  See [] for details.

Our main example has $\theta' = \frac{1}{6}$ and $\theta'' = \frac{6}{7}$, so $\mathcal{J}'$ is the Dendrite discussed in detail in [ADS] and $\mathcal{J}''$ is the anti-Rabbit (the mirror image of the Rabbit $p=2, k=3$).  In this example $X'^{(m)}$ and $1 - X''^{(m)}$ consist of exactly the same points, only with different identifications.  Thus $X^{(m)}$ consists of all points of the form $\frac{2^n + 7 \ell}{7(2^m)}$ for $n = 1, 2, 3$ and $0 \leq \ell \leq 2^m - 1$.

The rabbit identification simply identifies $\frac{2+7\ell}{7(2^m)}$, $\frac{4+7\ell}{7(2^m)}$ and $\frac{8+7\ell}{7(2^m)}$.  We will call these triplets.  The Dendrite makes many of the same identifications, but also makes some different identifications that result in pairs of triplets begin identified to create {\it sextuplets}.   The first example occurs when $m=4$, with the identification of the $\ell = 1$ and $\ell = 9$ triplet.  Therefore, for $m \geq 4$, there will be identifications of pairs of triplets for all $\ell \equiv 1 \mod 2^3$.  For $m=6$ a different type of sextuple arises by identifying the $\ell = 5$ and $\ell = 7$ triplet, leading to identifications for all $\ell \equiv 5 \mod 2^5$ for $m \geq 6$.  Note that the triplets first identified for $m=4$ and $m=6$ contain the dividing points $\frac{1}{12}$ and $\frac{7}{12}$ for the dendrite.  More generally, for $m=2j+2$ there will be triplets corresponding to $\ell = \frac{1}{3}(2^{2j}-1)$ and $\ell = \frac{1}{3}(2^{2j}-1)+2^{2j+1}$ that contain these dividing points and so get identified in a sextuplet, and for all $m \geq 2j+2$ this will result in identifications of all triplets with $\ell \equiv \frac{1}{3}(2^{2j}-1) \mod 2^{2j+1}$.  The exact description of the pairing of such triplets is rather complicated and will not be given here.  These give all the sextuplets, and if we denote by $a_m$ the number on level $m$, then
$$
a_{m+1} =
\begin{cases} 2a_m & \text{if $m$ is even,}
\\
2a_m+1 &\text{ if $m$ is odd.}
\end{cases}
$$
which leads to,
$$
\begin{cases} a_{2m} = \frac{4^{m-1}-1}{3} & \text{}
\\
a_{2m+1} = \frac{2 \cdot 4^{m-1}-2}{3} &\text{ }
\end{cases}
$$
for $m\geq 2$

The graph $X^{(m)}$ consists of intervals of lengths $\frac{2^j}{7(2^m)}$ for $j=1, 2, 3$, and we will call these intervals of type 1, 2 and 3.  When we pass to the next graph $X^{(m+1)}$, the intervals of type 1 and 2 do not subdivide and simply become intervals of type 2 and 3, respectively.  Intervals of type 3 subdivide into four intervals of types 1, 2, 3, 1 in that order.  Some intervals will be loops (endpoints identified), and those do not contribute to the energy.  Those that are not loops will always be paired, meaning that there will be two (and sometimes four) intervals with identified endpoints.  Figure 3.1 shows the subdivision of paired intervals in $X^{(4)}$ when increasing $m$ to 5, yielding two sets of paired type 1 intervals.  A more detailed analysis of identifications in subdivision will be given in section 5.  Here the important observations is that the intervals of type 2 and 3 created in the subdivision are loops and so do not contribute to energy, and intervals of type 1 are all paired.

\begin{figure}[h]
\numberwithin{figure}{section}
\caption{}
\centering
\includegraphics[scale = 0.6]{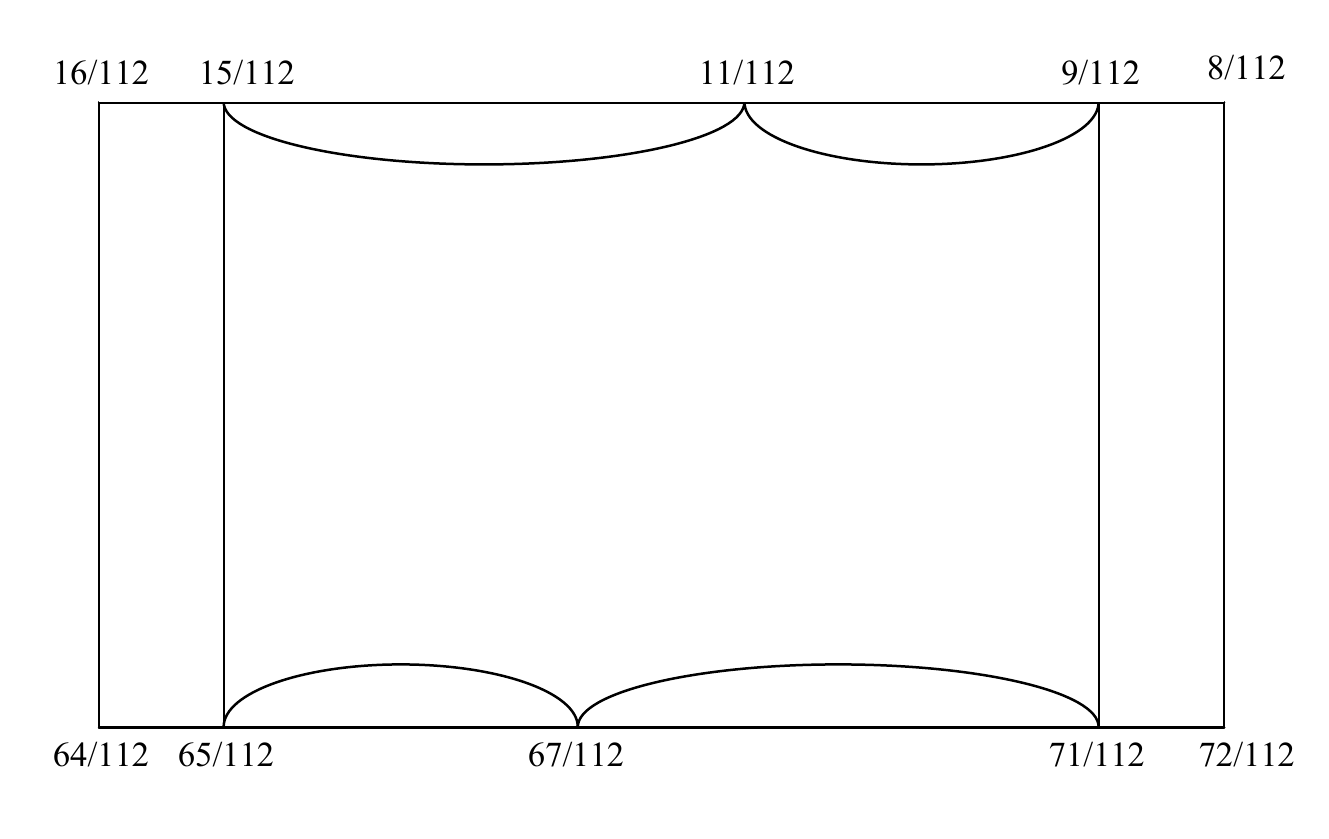}
\end{figure}

The analysis of the energy is then exactly the same as for the Rabbit.  We consider the three energies $E_1^{(m)}(u)$, $E_2^{(m)}(u)$ and $E_3^{(m)}(u)$ and observe that $E_3^{(m)}(\tilde{u}) = E_2^{(m-1)}(u)$, $E_2^{(m)}(\tilde{u}) = E_1^{(m-1)}(u)$ and $E_1^{(m)}(\tilde{u}) = 4E_3^{(m-1)}(u)$, so the matrix is
$$
M = \begin{pmatrix}
0 & 0 & 4 \\ 
1 & 0 & 0 \\
 0 & 1 & 0
 \end{pmatrix}
$$

Thus,

 \begin{equation*}
\mathcal{E}^{(m}(u) = 4^{-\frac{m}{3}}(E_1^{(m)}(u)+4^{\frac{1}{3}}E_2^{(m)}(u)+4^{\frac{2}{3}}E_3^{(m)}(u)).
\tag{3.1}\end{equation*} 

Perhaps the simplest nontrivial matings are the Basilica/Rabbit and Rabbit/Rabbit matings.  We studied these examples extensively searching for a way to sort the edges into a finite number of types, and then find subdivision rules for passing from one graph to the next.  What we found was that each time we increased the level $m$, it was necessary to add new edge types.  It seems unlikely that there would be only a finite number of edge types that would work for all $m$, although there is the remote possibility that we stopped the search too soon.  We then investigated constructing energies using a small number of edge types and weights chosen through trial and error, attempting to obtain approximate Laplacians that had some rough consistency from level to level.  The results were not very satisfactory.  Some data from these attempts may be found on the website.

\pagebreak

\section{Numerical Data}

In this section we present numerical data on the examples above.  In each case, to obtain the data, we approximated the Julia set, or formal mating of Julia sets, by a graph.  We then computed the graph Laplacian on the approximating graph, with the adjustment to the standard graph Laplacian coming from the weighting of the edges of the graph, and the measure placed on the graph.  Using the built in eigenvalue and eigenvector functionalities in Mathematica and MATLAB, we produced approximations of the eigenvalues and eigenfunctions of the Laplacian on the set.

The data presented comes from the highest level approximation computed.  We present a selection of eigenfunctions, the first 150 eigenvalues, the eigenvalue counting function and the Weyl ratio.

The first example, the cubic Basilica, corresponds to $p=3$ and $k=2$, so $\theta = \frac{1}{24}$, the first Julia set shown in figures 2.1 and 2.2.  Figure 4.1 shows the first eight non-constant eigenfunctions, figure 4.2 shows the eigenvalue counting function, $N(t) = \#\{\lambda_j : \lambda_j \leq t\}$, and figure 4.3 shows the Weyl ratio $W(t) = N(t)/t^{\alpha}$ for the appropriate power $\alpha$, which in this case is $\alpha = \frac{2}{3}$.  Table 4.1 shows the first 150 eigenvalues.  The second example, the cubic Rabbit, corresponds to $p=3$ and $k=3$, so $\theta = \frac{1}{78}$, the second Julia set show in figures 2.1 and 2.2.  The corresponding data here is shown in figures 4.4-6 (with $\alpha = \frac{3}{4}$) and table 4.2.

The last example is the dendrite-antiRabbit formal mating discussed in section 3.  Figure 4.7 shows the first eight eigenfunctions and figure 4.8 shows a multiplicity 5 eigenspace.  Figure 4.9 shows the eigenvalue counting function and figure 4.10 shows the the Weyl ratio with the value $\alpha = \frac{7}{10}$ determined experimentally.  Table 4.3 shows the first 150 eigenvalues.

\pagebreak

\includepdf[pages=-]{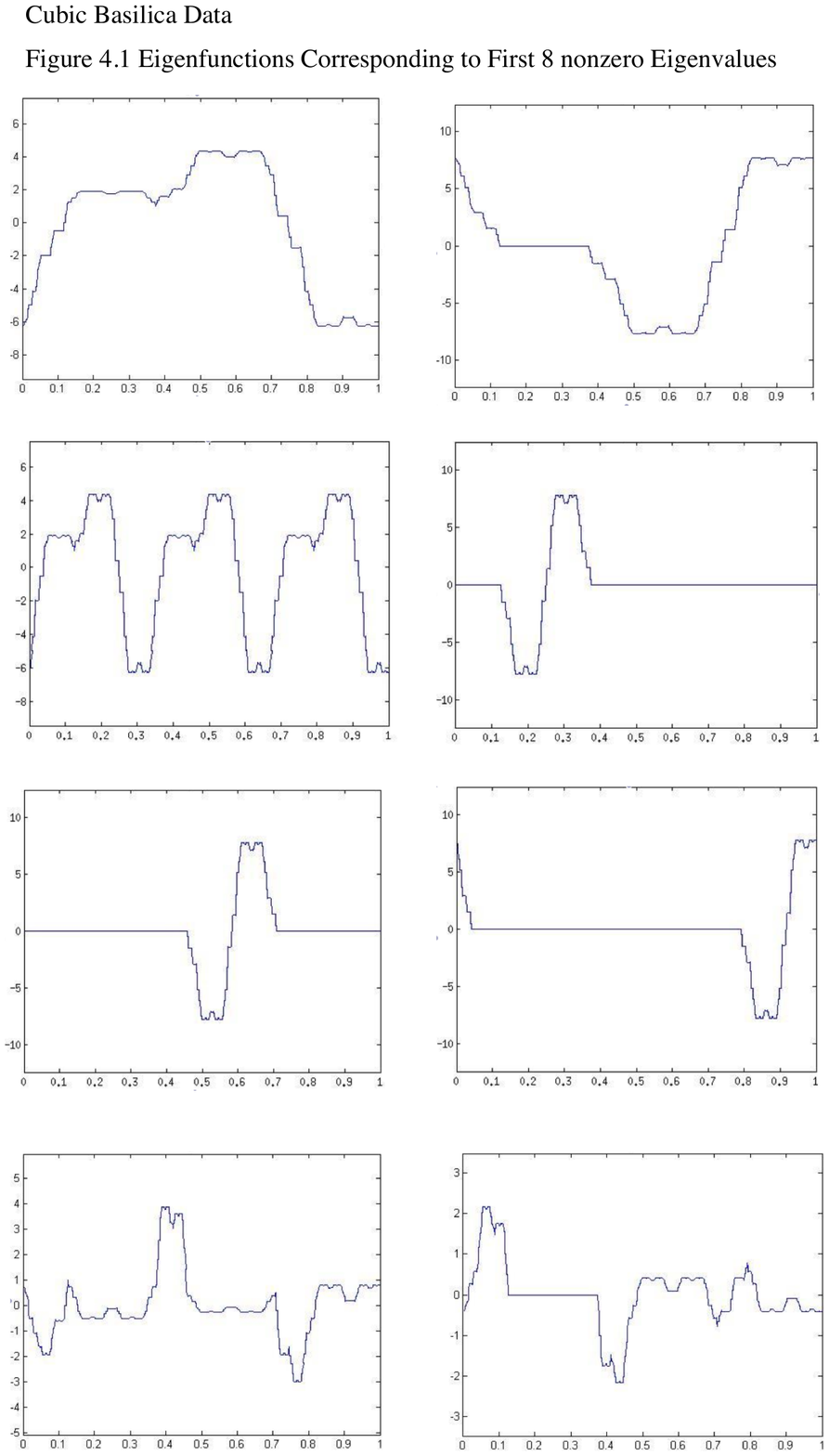}
\pagebreak

\section{Properties of the Spectrum}

We consider first the simple Julia sets discussed in section 2 with parameters $(p, k)$.  Note that (2.15) implies that if $u$ is an eigenfunction of $-\Delta$ with eigenvalue $\lambda$, then $u \circ P$ is an eigenfunction with eigenvalue $p^{\frac{k+1}{k}}\lambda$.  In particular, the spectrum is preserved under multiplication by $p^{\frac{k+1}{k}}$.  It therefore makes sense to split the nonzero eigenvalues into the {\it primitive} ones, not of the form $p^{\frac{k+1}{k}}\lambda'$ for some other eigenvalue $\lambda'$, and the {\it derived} ones that have this form.  Moreover, the derived eigenvalues have an{\it order} $j$, the unique value for which the eigenvalue for which the eigenvalue $p^{j\frac{k+1}{k}}\lambda'$ for $\lambda'$ a primitive eigenvalue.  Note that primitive eigenvalues will have order zero under this definition.  Let $m(\lambda)$ denote the multiplicity of the $\lambda$-eigenspace

\newtheorem{thm}{Conjecture}[section]
\begin{thm}
 Let $p=3$.  Then $m(\lambda)$ only depends on the order of $\lambda$, $m(\lambda) = m_j$ where $m_0 = 2$ and the recursion relations

\begin{equation*}
m_j =
\begin{cases} 3(m_{j-1}-1)+1 & \text{if $j \equiv 1, 2, ..., k-1 \mod k$,}
\\
3(m_{j-1}-2)+2 &\text{ if $j \equiv 0 \mod k$.}
\end{cases}
\tag{5.1}\end{equation*}

hold.  In particular, all multiplicities are even (for $\lambda \neq 0$).
 \end{thm}
 
 We will sketch an argument for why we believe the conjecture is valid.  It is unlikely that this can be made into a proof since it uses the "no coincidences" principle: unrelated eigenfunctions will have different eigenvalues.  (A similar principle may be invoked in describing the spectrum of the ordinary Laplacian on the unit disc, where it amounts to saying that Bessel functions of distinct orders do not have common zeros.)   The occurrence of coincidences would increase the multiplicities (but they would still be even).  Because of the spectral clustering phenomena discussed below, we expect that there will be many near coincidences, so numerical approximations will not be effective in ruling out coincidences.
 
 We note that the conjecture gives a description of multiplicities that differs substantially from the examples with $p=2$ discussed in [FS].  We expect that suitable modifications of the conjecture will be valid for odd values of $p$.
 
 The key observation is that there is a symmetry group of dihedral-3 type acting on these Julia sets.  It is clear that the rotation $t \rightarrow t+\frac{1}{3}$ is a symmetry, and this generates a $\mathbb{Z}_3$ symmetry group.  When $k = 2$ it is also clear that the reflections of the circle about the diameters $(\frac{1}{4}, \frac{3}{4})$, $(\frac{1}{12}, \frac{7}{12})$, $(\frac{5}{12}, \frac{11}{12})$ are symmetries of $\mathcal{J}$.  For $k \geq 3$ there are "reflection" symmetries that are not true reflections of the circle, but still fix one of the points in $X^{(1)}$ and permute the other two.  For example, this map will actually be the identity on the whole arc containing the identified points on the circle corresponding to the fixed points in $X^{(1)}$, in other words on the set of "ears" of $\mathcal{J}$ that attach to the point.  A similar idea was discussed for $p=2$ in [FS].
 
 Symmetries preserve the eigenspaces, so each eigenspace splits into a direct sum of spaces that transform according to the irreducible representations of the symmetry group.  In this case there are just three such representations: two one-dimensional representations, the trivial and alternating representations, and one two-dimensional representation.  Note that both one-dimensional representations act like the identity on the rotational $\mathbb{Z}_3$ subgroup.
 
 \newtheorem{thm2}[thm]{Lemma}
 \begin{thm2}
 An eigenfunction $u$ is of the form $v \circ P$ for another eigenfunction $v$ if and only if 
 
 \begin{equation*}
u(t+\frac{1}{3}) = u(t)
\tag{5.2}\end{equation*}
 
 As a consequence, a primitive eigenspace cannot contain any components corresponding to one-dimensional representations.
 \end{thm2}
 
\begin{proof}
If (5.2) holds, we may define $v(t) = u(\frac{1}{3}t)$.  Then $u(t) = v(3t)$ and $v$ is an eigenfunction.  Conversely, if $u(t) = v(3t)$, then clearly (5.2) holds.
\end{proof}
 
 This shows that primitive eigenspaces have even multiplicity, and the no coincidences principle would make $m_0 = 2$.  However, we can arrive at another important insight with the observation that the two dimensional representation contains one dimensional subspaces that are skew-symmetric under one of the reflections.  Thus, the eigenspace contains a function $u$ that satisfies $u \circ R = -u$ for the reflection that fixes a point $x$ in $X^{(1)}$.  Now $x$ splits $\mathcal{J}$ into a large arc and a small arc, and because of the skew-symmetry we can write $u$ as a sum of eigenfunctions supported on each of them.  However, we may argue that the one supported on the small arc must be zero, for if not we could rotate it around to create a nonzero eigenfunction satisfying (5.2), contradicting the fact that we have a primitive eigenspace.  Therefore we conclude that the primitive eigenspace has a function supported on the large arc and vanishing on the small arc.  Taking two points in $X^{(1)}$ gives us a basis (of course the functions associated with all three points are linearly dependent).  We already observe a contrast with the $p=2$ case discussed in [FS] where essentially every eigenfunction vanishes on an open set, since if we take a sum of two such functions the result will not have such a vanishing property.
 
 Next we describe how to pass from a primitive eigenspace to the derived eigenspace of order 1.  Let $u$ be the eigenfunction that vanishes on $[\frac{1}{3^k-1}, \frac{3^{k-1}}{3^k-1}]$.  Then $v(x) = u(3x)$ vanishes on $[\frac{1}{3(3^k-1)}, \frac{3^{k-1}}{3(3^k-1)}]$, $[\frac{3^k }{3(3^k-1)}, \frac{3^k+ 3^{k-1}-1}{3(3^k-1)}]$ and $[\frac{2\cdot3^k-1 }{3(3^k-1)}, \frac{2\cdot3^k+ 3^{k-1}-2}{3(3^k-1)}]$.  Note that $v$ is supported on the union of three intervals $I = [\frac{3^{k-2}}{3^k-1}, \frac{3^{k-1}}{3^k-1}]$, $I+\frac{1}{3}$ and $I+\frac{2}{3}$, and these intervals are separated in $\mathcal{J}$.  Therefore we can split $v$ into a sum of three eigenfunctions supported in each of these intervals.  On the other hand $v(x + \frac{1}{9}) = u(3x+\frac{1}{3})$ is supported on $I+\frac{2}{9}$, $I+\frac{5}{9}$ and $I+\frac{8}{9}$ and these intervals are connected in $\mathcal{J}$.  Thus we arrive at the conclusion that $m_1 = 4$ with a basis of eigenfunctions vanishing on certain open intervals.
 
 The process of passing from derived eigenspaces of order $j-1$ to derived eigenspaces of order $j$ is similar.  If $(u_1, ..., u_{m_{j-1}})$ is our special basis for the order $j-1$ eigenspace, we consider the functions $v_m(x) = u_m(3x)$.  Depending on the residue of $j \mod k$, all but one or all but two of them will be supported on separated intervals, and may be split into three linearly independent eigenfunctions.  This gives the relation (5.1), and a special basis.
 
 We can make another conjecture concerning the location in the spectrum of derived eigenvalues.  Let $\lambda_0 = 0$ denote the lowest eigenvalue, and then label consecutive eigenvalues $\lambda_1, \lambda_2, \lambda_3, ...$ repeated according to multiplicity.  For any eigenvalue $\lambda$ denote by $\#(\lambda)$ the largest $n$ such that $\lambda_n = \lambda$.  We know that the eigenvalue counting function $N(t)$ has growth rate $t^\alpha$ for $\alpha = \frac{\log p^{\frac{k}{k+1}}}{\log p} = \frac{k}{k+1}$ and the Weyl ration $\frac{N(t)}{t^\alpha}$ tends to a multiplicatively periodic function of period $p^{\frac{k+1}{k}}$.  That means $\#(p^{j\frac{k+1}{k}}\lambda) \approx p^j \#(\lambda)$.
 
 \newtheorem{thm3}[thm]{Conjecture}
 \begin{thm3}
 Let $p=3$.  Then 
 
  \begin{equation*}
\#(3^{j\frac{k+1}{k}}\lambda) = 3^j\#(\lambda).
\tag{5.3}\end{equation*}
 
 \end{thm3}
 
 This conjecture is based solely on experimental evidence.  It is not valid for $p=2$.  A similar conjecture was described in [ADS] for the Dendrite example, but not for the other examples presented there.  A similar result was proved for the Sierpinski gasket in [S].  The experimental evidence points to the existence of spectral gaps between $\lambda_{3^n-1}$ and $\lambda_{3^n}$, although it is not clear what the lower bound of the ratio $\frac{\lambda_{3^n-1}}{\lambda_{3^n}}$ might be.  Also, we see many eigenvalues that are distinct relatively but close together, suggesting that spectral clustering (arbitrarily many distinct eigenvalues in arbitrarily small intervals) may occur.
 
 Next we discuss the spectrum for the Dendrite-Anti Rabbit mating.  For this we need a precise definition of the two reflection symmetries.  For the rabbit there are reflections we called {\it horizontal} $\rho_H$ and {\it vertical} $\rho_V$ in [FS], so we will keep the same terminology here, since the reflection on the mating are essentially the same.  (Note that these are not symmetries of the Dendrite.)  The horizontal reflection is the map $t \rightarrow t+\frac{1}{2}$ on the intervals $[\frac{1}{7}, \frac{4}{7}]$ and $[\frac{9}{14}, \frac{1}{14}]$ and maps $[\frac{1}{14}, \frac{1}{7}]$ to itself, and also $[\frac{4}{7}, \frac{9}{14}]$ to itself, while the vertical reflection is the identity on the intervals $[\frac{1}{7}, \frac{4}{7}]$ and $[\frac{9}{14}, \frac{1}{14}]$ and interchanges the intervals $[\frac{1}{14}, \frac{1}{7}]$ and $[\frac{4}{7}, \frac{9}{14}]$.

\begin{figure}[h]
\numberwithin{figure}{section}
\caption{}
\centering
 \includegraphics[scale = 0.75]{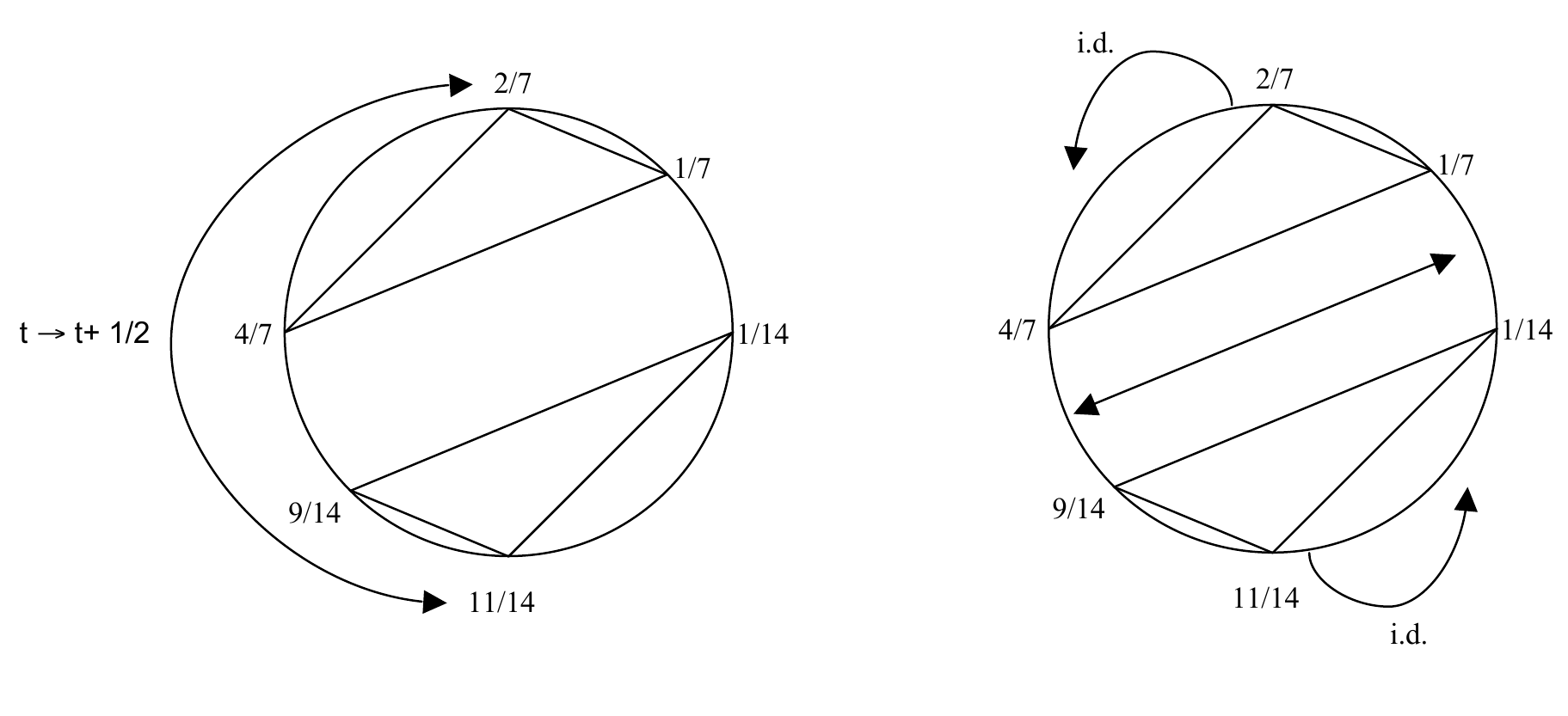}
\end{figure}

To describe those reflections on the remaining intervals $[\frac{1}{14}, \frac{1}{7}]$ and $[\frac{4}{7}, \frac{9}{14}]$ we observe that we may define an "orientation reversing" map between two intervals in the $X^{(m)}$ graph of type 3 iteratively as follows.  Say the intervals are $[\frac{4+7\ell}{7(2^m)}, \frac{8+7\ell}{7(2^m)}]$ and $[\frac{4+7\ell'}{7(2^m)}, \frac{8+7\ell'}{7(2^m)}]$, and we want to map $\frac{4+7\ell}{7(2^m)}$ to $\frac{8+7\ell'}{7(2^m)}$ and $\frac{8+7\ell}{7(2^m)}$ to $\frac{4+7\ell'}{7(2^m)}$.  We subdivide the intervals in $X^{(m+1)}$ 
\begin{figure}[h]
\numberwithin{figure}{section}
\caption{}
\centering
\includegraphics[scale = 0.75]{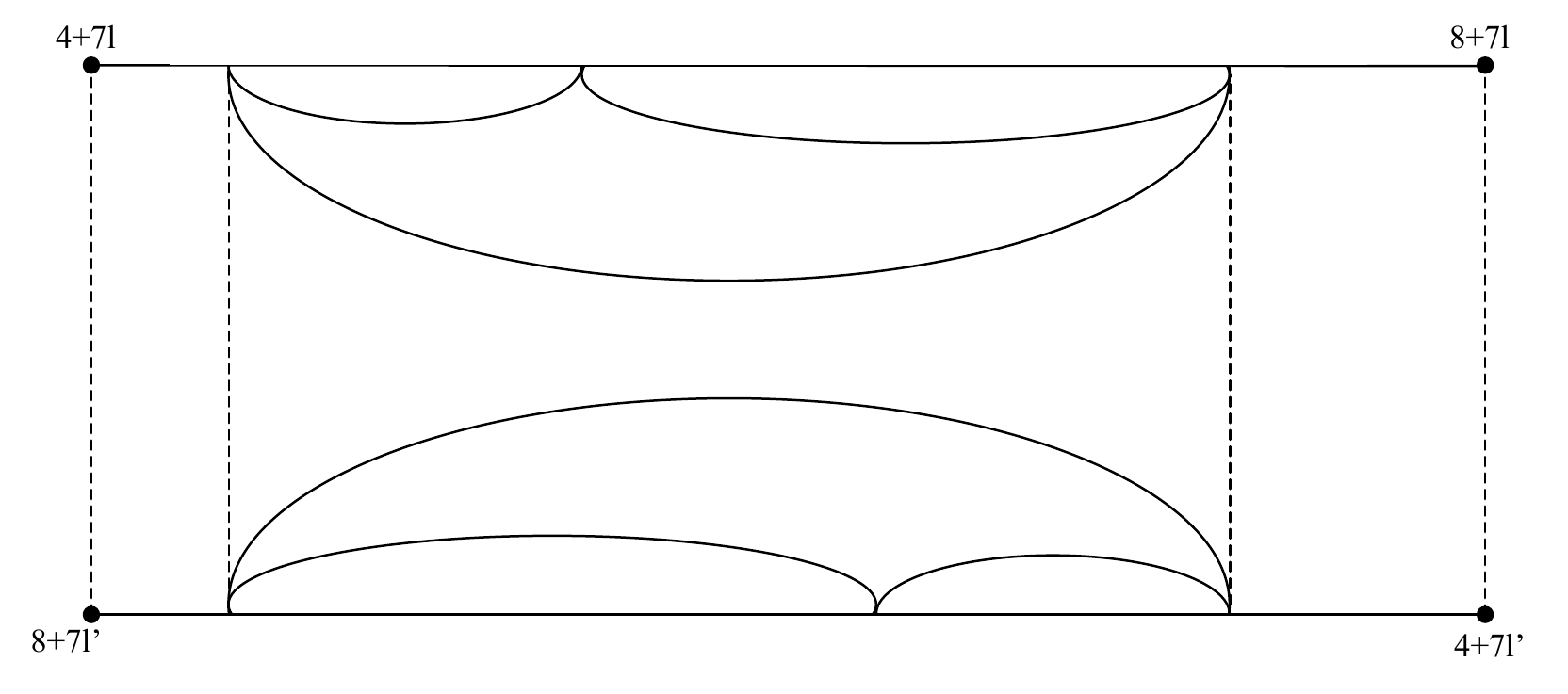}
\end{figure}
and map the outer subintervals to each other in orientation reversing fashion, and the inner subintervals are mapped to each other by translation.  See Figure 5.2.  In particular, $\rho_H$ is the identity on $[\frac{9}{112}, \frac{15}{112}]$ and maps $[\frac{8}{112}, \frac{9}{112}]$ to $[\frac{15}{112}, \frac{16}{112}]$ reversing orientation (and similarly on $[\frac{4}{7}, \frac{9}{14}]$), and $\rho_V$ is the translation $t \rightarrow t+\frac{1}{2}$ mapping $[\frac{9}{112}, \frac{15}{112}]$ to $[\frac{65}{112}, \frac{71}{112}]$ (and backward), while it is orientation reversing between $[\frac{8}{112}, \frac{9}{112}]$ and $[\frac{71}{112}, \frac{72}{112}]$ and similarly between $[\frac{15}{112}, \frac{16}{112}]$ and $[\frac{64}{112}, \frac{65}{112}]$.
 
 These two reflections generate a symmetry group isomorphic to $\mathbb{Z}_2 \times \mathbb{Z}_2$, the composition in either order being the translation $t \rightarrow t+\frac{1}{2}$.  Primitive eigenfunctions must be skew-symmetric with respect to this translation, so they fall into two types: {\it horizontal} eigenfunctions are skew-symmetric with respect to $\rho_H$ and symmetric with respect to $\rho_V$, while {\it vertical} eigenfunctions have the reverse symmetry.  Note that horizontal eigenfunctions vanish on $[\frac{9}{112}, \frac{15}{112}]$ and $[\frac{65}{112}, \frac{71}{112}]$ because $\rho_H$ is the identity on those intervals.  It then follows that eigenfunctions derived from horizontal eigenfunctions will vanish on prescribed intervals, but will generate one-dimensional eigenspaces (assuming no coincidences).
 
 Similar reasoning shows that vertical eigenfunctions are supported in $[\frac{1}{14}, \frac{1}{7}] \cup [\frac{4}{7}, \frac{9}{14}]$, but actually we can say more.
 
  \newtheorem{thm4}[thm]{Theorem}
 \begin{thm4}
Primitive vertical eigenfunctions are supported in $[\frac{11}{112}, \frac{15}{112}] \cup [\frac{67}{112}, \frac{71}{112}]$
 \end{thm4}
 
 \begin{proof}
 If $\rho_V(t) \sim t$ then $u(\rho_V(t)) = -u(t)$ imples $u(t) = 0$.  From this we infer that $u$ is supported on $[\frac{1}{14}, \frac{1}{7}] \cup [\frac{4}{7}, \frac{9}{14}]$, and it also vanishes at the sextuplet $\frac{9}{112} \sim \frac{11}{112} \sim \frac{15}{112} \sim \frac{65}{112} \sim \frac{67}{112} \sim \frac{71}{112}$.  It also follows from the skew-symmetry that we may cut the eigenfunciton at these points, to obtain three distinct eigenfunctions supported on $[\frac{8}{112}, \frac{9}{112}] \cup [\frac{15}{112}, \frac{16}{112}] \cup [\frac{64}{112}, \frac{65}{112}] \cup [\frac{71}{112}, \frac{72}{112}]$, on $[\frac{9}{112}, \frac{11}{112}] \cup [\frac{65}{112}, \frac{67}{112}]$, and on $[\frac{11}{112}, \frac{15}{112}] \cup [\frac{67}{112}, \frac{71}{112}]$.  We need to show that the first two vanish identically, and for this we need to use the fact that we are assuming the eigenvalue is primitive.  Suppose first that $u$ is supported on $[\frac{9}{112}, \frac{11}{112}] \cup [\frac{65}{112}, \frac{67}{112}]$.  Then
$$
v(x) =
\begin{cases} u(\frac{1}{4}x-\frac{1}{16}) & \text{if $x \in [\frac{4}{7}, \frac{9}{14}]$,}
\\
-u(\frac{1}{4}x+\frac{1}{16}) &\text{ if $x \in [\frac{1}{14}, \frac{1}{7}]$,}
\\
0 &\text{otherwise.}
\end{cases}
$$
is also an eigenfunction, contradicting the assumption that the eigenvalue is primitive.  Similarly, if $u$ is supported in the first union of four intervals,
$$
v(x) =
\begin{cases} u(\frac{1}{8}x) & \text{if $x \in [\frac{4}{7}, \frac{9}{14}]$,}
\\
-u(\frac{1}{8}x+\frac{1}{16}) &\text{ if $x \in [\frac{1}{14}, \frac{1}{7}]$,}
\\
0 &\text{otherwise.}
\end{cases}
$$
is also an eigenfunction.
\end{proof}
 
 We can now describe the derived eigenspaces and their multiplicities.  If $u(t)$ is a primitive vertical eigenfunction, then $u(2t)$ is a derived eigenfunction supported on $[\frac{11}{7(2^5)}, \frac{15}{7(2^5)}] \cup [\frac{67}{7(2^5)}, \frac{71}{7(2^5)}] \cup [\frac{123}{7(2^5)}, \frac{127}{7(2^5)}] \cup [\frac{189}{7(2^5)}, \frac{193}{7(2^5)}]$.  This support is not connected, however, and splits into the connected sets $[\frac{11}{7(2^5)}, \frac{15}{7(2^5)}]  \cup [\frac{189}{7(2^5)}, \frac{193}{7(2^5)}]$ and $[\frac{67}{7(2^5)}, \frac{71}{7(2^5)}] \cup [\frac{123}{7(2^5)}, \frac{127}{7(2^5)}]$.  By restricting $u(2t)$ to each of these sets we obtain a 2-dimensional eigenspace.  By iterating this argument we obtain a splitting of $u(2^kt)$ into $2^k$ linearly independent eigenfunctions.  But this is not the whole story.  There are two other types of eigenfunctions that arise.
 
 To describe the first type consider the case $k=2$.  Each of the four pieces of $u(4t)$ is supported in a pair of intervals of length $\frac{1}{7(2^4)}$, and the eight total intervals are obtained from $[\frac{11}{7(2^6)}, \frac{15}{7(2^6)}]$ by rotation by multiples of $\frac{1}{8}$.  However, there is another pair of intervals $[\frac{39}{7(2^6)}, \frac{43}{7(2^6)}]$ and $[\frac{263}{7(2^6)}, \frac{267}{7(2^6)}]$ that are paired and support another eigenfunction obtained from pieces of $u(4t)$ that are rotated by multiples of $\frac{1}{16}$.  Thus we get a multiplicity 5 eigenspace.  When $k = 3$ we double the multiplicity to 10 by splitting each of the 5 eigenfunctions from $k=2$.  For $k=4$ we split the 10 previous eigenfunctions, add one more by paired intervals rotated by multiples of $\frac{1}{2^5}$, but also there is yet another type of eigenfunction supported on a set of four intervals with endpoints belonging to two sextuplets of identified points $[\frac{36}{7(2^6)}, \frac{37}{7(2^6)}] \cup [\frac{43}{7(2^6)}, \frac{44}{7(2^6)}] \cup [\frac{260}{7(2^6)}, \frac{261}{7(2^6)}] \cup [\frac{267}{7(2^6)}, \frac{268}{7(2^6)}]$.  This is illustrated in Figure 5.3.  This is the $j=2$ case of a more general construction that yields an eigenfunction supported on  $2^j$ intervals with endpoints belonging to $2^{j-1}$ sextuplets of identified points.  These eigenfunctions arise at level $k$ for $k \geq 4$ even with $j \leq \frac{k}{2}$.

\begin{figure}[h]
\numberwithin{figure}{section}
\caption{}
\centering
 \includegraphics[scale = 0.75]{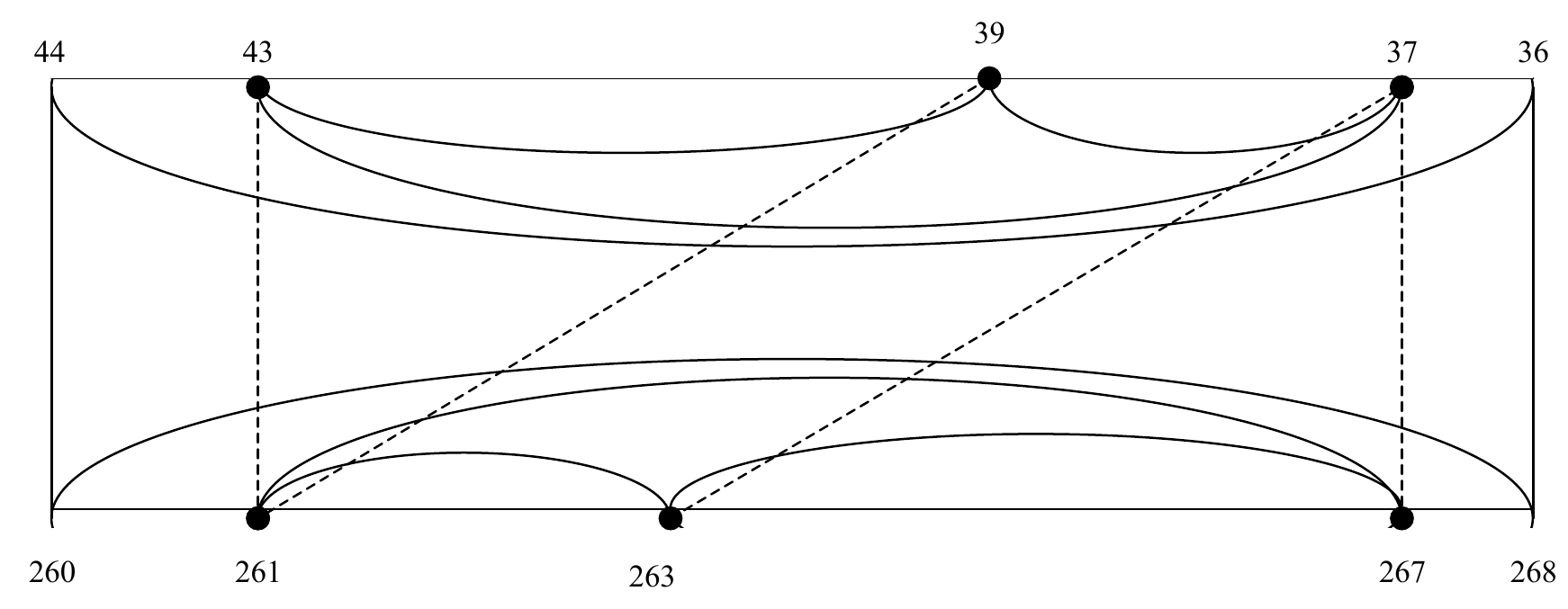}
\end{figure}

 In general, when $k$ is odd we just split the eigenfunctions from level $k-1$.  When $k$ is even (and $k \geq 4$) we also add a total of $\frac{k}{2}$ eigenfunctions are described above.  If $b_k$ denotes the multiplicities then
 
 $$
b_{k} =
\begin{cases} 2b_{k-1} & \text{if $k$ is odd,}
\\
2b_{k-1}+\frac{k}{2} &\text{ if $k$ is even, and $k \geq 4$,}
\end{cases}
$$

This implies $b_{2k} = 4b_{2k-2}+k$, which is easily solved to obtain $b_{2k} = \frac{13}{9}4^k - \frac{4}{9} - \frac{k}{3}$, and then $b_{2k+1} = 2(\frac{13}{9}4^k -\frac{4}{9} - \frac{k}{3})$.

\section{Acknowledgements}

We are grateful to John Hubbard and Dierk Schleicher for helpful discussions

\bibliographystyle{amsplain}

\end{document}